\documentclass[a4paper]{article}

\usepackage{color}
\usepackage{graphicx}
\usepackage{amsmath,amsfonts,amsthm,amssymb}
\usepackage{wrapfig}
\usepackage[english]{babel}
\usepackage{mathrsfs}
\usepackage{dsfont}
\usepackage{hyperref}

\newtheorem{lem}{Lemma}
\newtheorem{prop}{Proposition}
\newtheorem{rem}{Remark}

\newcounter{ris}
\renewcommand{\r}{\refstepcounter{ris}%
                  Fig. \arabic{ris}}

\begin{document}
\renewcommand{\proofname}{Proof.}

\title{On Point Coverings of Boxes in $\mathbb R^d$}

\author{A.\,V.~Akopyan}
\date{}

\maketitle

\begin{abstract}
Families of boxes in $\mathbb R^d$ are considered. In the paper an upper
bound on the size of a minimum transversal in terms of the space dimension and
the independence number of the given family was improved.
\end{abstract}

Given a family $\mathscr B$ of boxes with edges parallel to coordinate axes of
$\mathbb R^d$, let $\nu (\mathscr B)$ denote the maximal number of pairwise
disjoint members of $\mathscr B$ and $\tau(\mathscr B)$ denote the minimal
number of points in a set meeting all the members of $\mathscr B$. 

Hadwiger and Debruner in \cite{Had} raised problem of finding $f(n,d)$ ---
the supremum of $\tau(\mathscr B)$ taken over families $\mathscr B$ of
$d$-dimensional boxes with $\nu (\mathscr B)$.

It was observed in \cite{Had} that 
\begin{equation*}
 f(1,d)=1\text{ }\forall d\text{ and } f(n,1) = n \text{ }\forall n\text{.}
\end{equation*}

Also in \cite{Had} it was shown that
\begin{equation*}
f(n,2)\leqslant\frac{n(n-1)}{2}\text{.}
\end{equation*}

First nontrivial lower bound was proved by Gy\'afr\'as and Lehel
\cite{GyLeh}

\begin{equation*}
1.5n\leqslant f(n,2)\text{.}
\end{equation*}

K\'arolyi \cite{Kar} proved that for any fixed $d$

\begin{equation*}
f(n,d) \leqslant (1+o(1))n\log^{d-1}_2n\text{.}
\end{equation*}

Fon-Der-Flaass and Kostochka \cite{FK} specified this upper bound. 
They showed that for any fixed~$d$

\begin{equation*}
 f(n,d)\leqslant n \log^{d-1}_2n + d-\frac{1}{2}n\log^{d-2}_2n
=(1+o(1))n\log^{d-1}_2n\text{.}
\end{equation*}

Also they showed that for $d=2$

\begin{equation*}
f(n,2)\geqslant \left[\frac{5n}{3}\right]\text{.}
\end{equation*}

For $d>2$ they proved
\begin{equation*}
 f(2,3)=4\text{, } f(2,4)=5\text{,  and } f(2,d)\geqslant
\frac{c'd^\frac{1}{2}}{\log d}\text{ for some $c'$.}
\end{equation*}

Applying an idea from \cite{FK} we improve the upper bound by $\thickapprox
1.057$, i.e. we show that

\begin{equation*}
f(n,d) < (\log_{\sqrt[3]{9}}2+o(1))n\log^{d-1}_{2}n\text{.}
\end{equation*}

In particular for $d=2$ and $n\geqslant 5$

\begin{equation*}
f(n,2)\leqslant n\log_{\sqrt[3]{9}}n\text{.}
\end{equation*}

\bigskip

Let us put, for convenience,
\begin{equation*}
 f(0,d)=0 \text{.}
\end{equation*}

In  \cite{FK} there was such useful proposition.

\begin{prop}
  \label{prop:(flkost)}
  For any $n\geqslant 2$ and $d \geqslant 2$

        \begin{equation*}
 f(n,d)\leqslant \min_{0\leqslant k \leqslant n-2}
\{f(k,d)+f(n-k-1,d)\}+f(n,d-1)\text{.}
        \end{equation*}

\end{prop}

Now we prove one simple corollary from that proposition.
It is needed for the bound for $f(n,d)$ in case of $d\geqslant3$.

\begin{lem}
  \label{lem:(fromflaas)} 
For any $n\geqslant 1$ and $d \geqslant 2$ 

 \begin{equation*}
    f(n,d)\leqslant n+\log_2{n}\cdot f(n,d-1)\text{.}
 \end{equation*}

\end{lem}

\begin{proof}
For $n=1$ that statement is obvious.
Suppose that it holds for all $n'<n$.
Note that from the proposition~\ref{prop:(flkost)} and the base of induction
we have

\begin{multline*}
f(n,d)\leqslant
f\left( \left\lceil \frac{n-1}{2}
\right\rceil,d \right)+f\left( \left\lceil \frac{n-1}{2}
\right\rceil,d \right)+f(n,d-1)=\\
=f(n,d-1)+2f\left(\left[\frac{n}{2}\right],d\right)
\leqslant \\ \leqslant 
f(n,d-1)+2\log_2{\left[\frac{n}{2}\right]}
\cdot
f\left(\left[\frac{n}{2}\right],d-1\right)+2\left[\frac{n}{2}\right]\text{.}
\end{multline*}

It is clear that $f(l,d)+f(m,d)\leqslant
f(l+m,d)$.
It follows that 
\begin{equation*}
 2f\left(\left[\frac{n}{2}\right],d-1\right) \leqslant
f(n,d-1)\text{.}
\end{equation*}

Therefore,

\begin{equation*}
f(n,d)\leqslant
%f(n,d-1)+2\log_2{\left[\frac{n}{2}\right]}
%\cdot f\left(\left[\frac{n}{2}\right],d-1\right)+2\left[\frac{n}{2}\right]
%\leqslant
n+(\log_2\frac{n}{2}+1)f(n,d-1)=n+\log_2{n}\cdot f(n,d-1)\text{.}
\end{equation*}

\end{proof}
The next lemma is needed for the bound in case of $d=2$.

\begin{lem}
  \label{lem:f(n,2,2)}
Suppose that for the family $\mathscr B$ of boxes in $\mathbb R^2$ there 
exists two prarllel to the absciss lines $l_1$ and $l_2$, such that
every box from $\mathscr B$ intersect at least one of them. Then
$\tau (\mathscr B)\leqslant \left[ \frac{3\nu (\mathscr B)}{2} \right]$.
\end{lem}

\begin{proof}
The proof is by induction on $\nu (\mathscr B)$ with a step of induction equal
to $2$.
In \cite{Had} it was shown that $f(1,2)=1$ and $f(2,2)=3$. 
This is base of induction.
Suppose this statement is proved for all families $\mathscr B$ with $\nu
(\mathscr B) \leqslant n-2$.
Let us assume that $\mathscr B= \{\mathscr B_i | i\in I\}$ with $\nu
(\mathscr B) = n$ and $\tau (\mathscr B)=f(n,b)$.
We may assume

\begin{equation*}
\mathscr B_i=[l_{i,x};r_{i,x}]\times[l_{i,y};r_{i,y}]\text{.}
\end{equation*}

For any real $x$ let

\bigskip

\begin{tabular}{ll}
  $\mathscr B^-(x)=\{\mathscr B_i|r_{i,x}<x\}\text{,}$ &
  $\mathscr B^+(x)=\{\mathscr B_i|l_{i,x}>x\}\text{,}$ \\
  $\mathscr B^0(x)=\mathscr B \setminus(\mathscr
B^-(x)\cup\mathscr B^+(x))\text{.}$  
\end{tabular}
\bigskip

Let $a=\sup \{x\in \mathbb R|\nu(\mathscr B^-(x))\leqslant 1\}$. 
By the choice of $a$ there exists $\mathscr B_a \in \mathscr
B^0(a)$ such that
$r_{a,x}=a$ and a box from $\mathscr B^-(a)$ disjoint with the box $\mathscr
B_a$

Therefore $\nu (\mathscr B^+(a))\leqslant n-2$, and by 
inductive assumption $\tau (\mathscr B^+(a))\leqslant [\frac{3(n-2)}{2}]$. 
Since in $\nu (\mathscr B^-(a))$ all boxes meet one another, we have $\tau
(\mathscr B^-(a))=1$. 
As well $\tau (\mathscr B^0(a))\leqslant 2$ because every box from
$\mathscr B^0(a)$ contains at least one of the point of intersection  of the
line $x=a$ and the lines $l_1$ and $l_2$.

\begin{equation*}
\tau (\mathscr B) \leqslant \tau (\mathscr B^-(a))+ \tau (\mathscr B^+(a))+ \tau
(\mathscr B^0(a))  \leqslant 1+ \left[\frac{3(\nu (\mathscr B)-2)}{2}
\right]+2=\left[ \frac{3\nu(\mathscr B)}{2} \right] \text{.}
\end{equation*}
\end{proof}

\begin{rem}
The bound from lemma \ref{lem:f(n,2,2)} is sharp.
\end{rem}

\begin{proof}
We need to construct a family $\mathscr B$ such that
$\tau (\mathscr B)=[\frac{3n}{2}]$.
For $n=1$ that family $\mathscr B$ consists of one box meeting at
least one of two lines $l_1$ and $l_2$.
The example for $n=2$ is the family $\mathscr B$ which consists of five
boxes and
$\tau (\mathscr B)=3$ (fig.~\ref{ris:f(3,2)}).
\begin{center}
\includegraphics{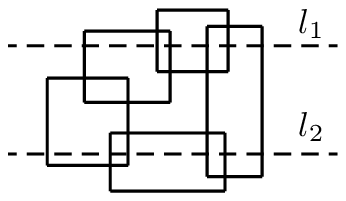}

\r\label{ris:f(3,2)}
\end{center}

For any $n$ the example consists of $[\frac{n}{2}]$ disjoint copies of the
family $\mathscr B^2$.
If $n$ is odd, then we can add one box disjoint with other
boxes.
 
\end{proof}

Applying lemma~\ref{lem:f(n,2,2)} we prove the next proposition in the same way
as proposition $\ref{prop:(flkost)}$. 

\begin{prop}
  \label{prop:(3chasti)}
  For any $n\geqslant 2$
 
\begin{equation*}
 f(n,2)\leqslant \min_{k+l+m=n-2}
\{f(k,2)+f(l,2)+f(m,2)\}+\left [\frac{3n}{2}\right]\text{.}
 \end{equation*}

\end{prop}

\begin{proof}

Consider $\mathscr B= \{\mathscr B_i | i\in I\}$ with $\nu (\mathscr B)=n$ and
$\tau(\mathscr B)=f(n,d)$. 
We may assume 
\begin{equation*}
\mathscr B_i=[l_{i,x};r_{i,x}]\times[l_{i,y};r_{i,y}]\text{.}
\end{equation*}

For any real $x_1, x_2 \in \mathbb R$ such that $x_1\leqslant x_2$ let
\bigskip

\begin{tabular}{l}
  $\mathscr B^-(x_1,x_2)=\{\mathscr B_i|r_{i,x}<x_1\}\text{,}$ \\
  $\mathscr B^\pm(x_1,x_2)=\{\mathscr B_i|l_{i,x}>x_1\text{
and }r_{i,x}<x_2\}\text{,}$\\
  $\mathscr B^+(x_1,x_2)=\{\mathscr B_i|l_{i,x}>x_2\}\text{,}$ \\
  $\mathscr B^0(x_1,x_2)=\mathscr B \setminus(\mathscr B^-(x_1,x_2)\cup
\mathscr B^\pm(x_1,x_2)\cup\mathscr B^+(x_1,x_2))\text{.}$
\end{tabular}

\bigskip

For any $k$, $l$ and $m$ such that $k+l+m= n-2$, denote
$a=\sup\{x\in \mathbb R|\nu(\mathscr B^-(x))\leqslant k\}$ and $b=\inf\{x\in
\mathbb R|\nu(\mathscr B^+(x))\leqslant m\}$. 
Note that $a\leqslant b$ otherwise  $\nu (\mathscr B)\leqslant k+m<n$.
It is clear that the family $\{\mathscr B_i | r_{i,x}\leqslant a\}$ consists of
$k+1$ disjoint boxes, therefore $\nu(\mathscr B^\pm(a,b) \cup
\mathscr B^+(a,b))\leqslant n-k-1=l+m+1$.
Similarly, $\{\mathscr B_i | l_{i,x}\geqslant b\}$ consists of at least $m+1$
disjoint boxes, therefore $\nu(\mathscr B^\pm(a,b))\leqslant n-k-1-m-1=~l$.
Thus $\tau(\mathscr B^-(a,b))\leqslant f(k,2)$,
$\tau(\mathscr B^\pm(a,b))\leqslant f(l,2)$,  $\tau(\mathscr B^+(a,b))\leqslant
f(m,2)$.
And applying lemma~\ref{lem:f(n,2,2)} we see that 
$\tau(\mathscr B^0(a,b))\leqslant [\frac{3n}{2}]$. 
Hense,

\begin{multline*}
\tau (\mathscr B) \leqslant \tau (\mathscr B^-(a,b))+\tau (\mathscr B^\pm(a,b))+
\tau (\mathscr B^+(a,b)) + \tau (\mathscr B^0(a,b))  \leqslant \\
\leqslant f(k,2)+f(l,2)+f(m,2)+\left[\frac{3n}{2}\right] \text{.}
\end{multline*}

\end{proof}
Applying this proposition to $n=5$, we have  $f(5,2)\leqslant 10$, which is
$1$ less than the bound obtained by means of~\ref{prop:(flkost)}. 
Now we show that $f(n,2)$ bound better than $n\log_2n$.

\begin{prop}
  \label{prop:(otsenkf2d)}
   
    Let $h(n)=n\log_{\sqrt[3]{9}}n+n$. Then

\begin{equation*}
f(n,2)\leqslant h(n)\text{.}
\end{equation*}

\end{prop}

\begin{proof}
Note that $h(1)=f(1,2)$ and $h(2)\thickapprox3.892789>3=f(2,2)$. 
Suppose that our proposition is proved for all $n'<n$.
By proposition~\ref{prop:(3chasti)} we have

\begin{equation*}
f(n,2)\leqslant
3f\left(\left\lceil\frac{n-2}{3}\right
\rceil,2\right)+\left[\frac{3n}{2}\right]\text{.}
\end{equation*}

Since $h(n)$ increases monotonically for $n> e^{-1}$
and $\left\lceil\frac{n-2}{3}\right \rceil\leqslant\frac{n}{3}$, we have:

\begin{equation*} 
3f\left(\left\lceil\frac{n-2}{3}\right
\rceil,2\right)+\left[\frac{3n}{2}\right]\leqslant
3h \left(\frac{n}{3}\right)+\frac{3n}{2}=
n\log_{\sqrt[3]{9}}\frac{n}{3}+3\frac{n}{2}+\frac{3n}{2}=h(n) \text{.}
\end{equation*}

\end{proof}

\begin{rem}
Note that for $n\geqslant 5$, $f(n,2) < n\log_{\sqrt[3]{9}}n$. 
In case of $n<15$ it can be shown by calculating $h(n)$ and $f(n,2)$ using
proposition~\ref{prop:(3chasti)}.
For bigger $n$ the proof is similar to the proof of
proposition~\ref{prop:(otsenkf2d)}.
\end{rem}

Using proposition \ref{prop:(otsenkf2d)} and lemma~\ref{lem:(fromflaas)} we
show following.

\begin{prop}
 For any fixed $d\geqslant2$

\begin{equation*}
f(n,d)<(\log_{\sqrt[3]{9}}2+o(1))n\log^{d-1}_2n\text{.}
\end{equation*}

\end{prop}

\begin{proof}
For $d=2$ the statement follows from proposition~\ref{prop:(otsenkf2d)}. Suppose
that it holds for all $d'<d$.
By lemma~\ref{lem:(fromflaas)}, we have

\begin{multline*}
f(n,d)\leqslant n+\log_2{n}\cdot f(n,d-1)\leqslant
n+\log_2{n}\cdot(\log_{\sqrt[3]{9}}2+o(1))n\log^{d-2}_2n=\\
=(\log_{\sqrt[3]{9}} 2+o(1))n\log^{d-1}_2n\text{.}
\end{multline*}

\end{proof}

\end{document}